\documentclass{article}
\usepackage{amssymb}
\newtheorem{theorem}{Theorem}
\newtheorem{corollary}{Corollary}

\newtheorem{lemma}{Lemma}

\newenvironment{pf}{{\noindent{\it Proof. }}}{\ \rule{2mm}{2.5mm}\medskip}
\newenvironment{pft}{{\noindent{\it Proof of Theorem }}}{\ \rule{2mm}
{2.5mm}\medskip}

\mathchardef\za="710B  
\mathchardef\zb="710C  
\mathchardef\zg="710D  
\mathchardef\zd="710E  
\mathchardef\zve="710F 
\mathchardef\zz="7110  
\mathchardef\zh="7111  
\mathchardef\zvy="7112 
\mathchardef\zi="7113  
\mathchardef\zk="7114  
\mathchardef\zl="7115  
\mathchardef\zm="7116  
\mathchardef\zn="7117  
\mathchardef\zx="7118  
\mathchardef\zp="7119  
\mathchardef\zr="711A  
\mathchardef\zs="711B  
\mathchardef\zt="711C  
\mathchardef\zu="711D  
\mathchardef\zvf="711E 
\mathchardef\zq="711F  
\mathchardef\zc="7120  
\mathchardef\zw="7121  
\mathchardef\ze="7122  
\mathchardef\zy="7123  
\mathchardef\zf="7124  
\mathchardef\zvr="7125 
\mathchardef\zvs="7126 
\mathchardef\zf="7127  
\mathchardef\zG="7000  
\mathchardef\zD="7001  
\mathchardef\zY="7002  
\mathchardef\zL="7003  
\mathchardef\zX="7004  
\mathchardef\zP="7005  
\mathchardef\zS="7006  
\mathchardef\zU="7007  
\mathchardef\zF="7008  
\mathchardef\zW="700A  

\newcommand{\be}{\begin{equation}}
\newcommand{\ee}{\end{equation}}
\newcommand{\ra}{\rightarrow}

\newcommand{\bea}{\begin{eqnarray}}
\newcommand{\eea}{\end{eqnarray}}
\newcommand{\beas}{\begin{eqnarray*}}
\newcommand{\eeas}{\end{eqnarray*}}
\newcommand{\Z}{{\mathbb Z}}
\newcommand{\R}{{\mathbb R}}

\newcommand{\1}{{\mathbf 1}}

\newcommand{\D}{{\rm d}}
\newcommand{\Df}{{\rm d}^\zF}

\newcommand{\we}{\wedge}
\newcommand{\nn}{\nonumber}
\newcommand{\ot}{\otimes}

\newcommand{\pa}{\partial}
\newcommand{\ti}{\times}
\newcommand{\A}{{\cal A}}
\newcommand{\Li}{{\cal L}}

\newcommand{\Ll}{{\pounds}}
\def\lna{\lbrack\! \lbrack}
\def\rna{\rbrack\! \rbrack}
\def\rnaf{\rbrack\! \rbrack_\zF}
\def\rnah{\rbrack\! \rbrack\,\hat{}}

\def\zT{{\cal T}}

\begin{document}
\title{On characterization of Poisson and Jacobi structures
\footnote{Supported by KBN, grant No 2 P03A 041 18.}}
        \author{
        Janusz Grabowski \\
        Institute of Mathematics \\
                Polish Academy of Sciences \\
                ul. \'Sniadeckich 8, P.O.Box 137, 00-950 Warszawa, Poland\\
                {\tt jagrab@mimuw.edu.pl} \\
                        \\
                Pawe\l\ Urba\'nski \\
                Division of Mathematical Methods in Physics \\
                University of Warsaw \\
                Ho\.za 74, 00-682 Warszawa, Poland \\
                {\tt urbanski@fuw.edu.pl}}
\date{}
\maketitle
\begin{abstract} We characterize Poisson and Jacobi structures by  means
of complete lifts of the corresponding tensors: the lifts  have  to  be
related  to canonical structures by morphisms  of  corresponding  vector
bundles. Similar results hold for generalized Poisson and Jacobi
structures (canonical structures) associated  with  Lie algebroids and
Jacobi algebroids.


\end{abstract}

\section{Introduction}
Jacobi   brackets   are  local   Lie   brackets   on   the
algebra $C^\infty(M)$ of smooth functions on a manifold $M$.
This    goes    back    to    the  well-known  observation  by  Kirillov
\cite{Ki} that in the case  of $\A=C^\infty(M)$ every local Lie  bracket
on $\A$ is of first order (an algebraic version of  this fact  for
arbitrary commutative  associative algebra $\A$ has been proved in \cite{Gr}).

Since every  skew-symmetric first-order bidifferential operator  $J$  on
$C^\infty(M)$ is of the form $J=\zL+I\we\zG$, where $\zL$ is a bivector
field, $\zG$ is a  vector  field  and  $I$  is
identity, the corresponding bracket of functions reads
\be\label{jac}
\{ f,g\}_J=\zL(f,g)+f\zG(g)-g\zG(f).
\ee
The Jacobi identity for this bracket is usually written in terms of  the
Schouten-Nijenhuis bracket $\lna\cdot,\cdot\rna$ as follows:
\be\label{kk}
\lna\zG,\zL\rna=0,\quad
\lna\zL,\zL\rna=-2\zG\we\zL.
\ee
Hence, every Jacobi bracket on $C^\infty(M)$ can be identified with  the
pair $J=(\zL,\zG)$ satisfying the above conditions,  i.e.  with  a  {\it
Jacobi structure} on $M$ (cf. \cite{Li}).
Note that we use the version of  the  Schouten-Nijenhuis  bracket  which
gives  a  graded  Lie  algebra  structure  on  multivector  fields
and which differs from the classical one by signs.
The Jacobi bracket (\ref{jac}) has he following properties:
\begin{enumerate}
\item $\{ a,b\}=-\{ b,a\}$ (anticommutativity),
\item $\{ a,bc\}=\{   a,b\} c+b\{   a,c\}-\{a,\1\}bc$  (generalized
Leibniz rule),
\item $\{\{ a,b\},c\}=\{ a,\{   b,c\}\}-\{b,\{  a,c\}\}$ (Jacobi
identity),
\end{enumerate}

\noindent
The generalized Leibniz rule tells that the bracket is a  bidifferential
operator on  $C^\infty(M)$ of  first  order.  In  the  case  when
$\zG=0$  (or, equivalently, when the constant function $1$ is a central
element), we deal  with
a  {\it  Poisson bracket}  associated  with   the   bivector   field
$\zL$   satisfying $\lna\zL,\zL\rna=0$.

For a smooth manifold $M$ we denote by $\zL_M$ the canonical Poisson
tensor on $T^*M$, which in local Darboux coordinates $(x^l,p_j)$ has  the
form $\zL_M=\pa_{p_j}\we\pa_{x^j}$.
In \cite{GU} the  following  characterization  of  Poisson  tensors,  in
terms of the {\it complete (tangent) lift} of contravariant tensors
$X\mapsto  X^c$ from the manifold $M$ to $TM$, is proved.

\begin{theorem} A bivector field $\zL$ on a manifold $M$ is  Poisson  if
and only if  the  tensors  $\zL_M$  and  $-\zL^c$  on  $T^*M$  and  $TM$,
respectively, are $\sharp_\zL$-related, where
$$ \sharp_\zL:T^*M\ra TM,\quad \sharp_\zL(\zw_x)=i_{\zw_x}\zL(x).
$$
\end{theorem}
\noindent
So, for a Poisson tensor $\zL$ the map $\sharp_\zL:(T^*M,\zL_M)\ra
(TM,-\zL^c)$ is a Poisson map.

The aim of  this  note  is  to  generalize  the  above  characterization
including Jacobi brackets and {\it canonical structures} associated with
{\it Lie algebroids} and {\it Jacobi algebroids}.

\section{Lie and Jacobi algebroids}

A {\em Lie algebroid} is a vector bundle $\zt:E\ra M$, together with a
bracket $\lna\cdot,\cdot\rna$ on the $C^\infty(M)$-module  $Sec(E)$ of smooth
sections of $E$, and a bundle morphism $\zr :E\ra TM$ over the identity
on $M$, called the {\em  anchor} of the
Lie algebroid, such that \noindent
\begin{description}
\item{(i)} the bracket $\lna\cdot,\cdot\rna$ is $\R$-bilinear,
alternating, and satisfies the Jacobi identity;
\item{(ii)} $\lna X,fY\rna=f\lna X,Y\rna+\zr(X)(f)Y$ for all  $X,Y\in Sec(E)$
and  all $f\in C^\infty(M)$.
\end{description}
From (i) and (ii) it follows easily
\begin{description}
\item{(iii)} $\zr(\lna X,Y\rna)=[\zr(X),\zr(Y)]$ for all $X,Y\in Sec(E)$.
\end{description}
We will often identify sections $\zm$ of the dual bundle $E^*$ with
linear (along fibres)  functions  $\zi_\zm$  on the vector bundle $E$:
$\zi_\zm(X_p)=<\zm(p),X_p>$. If $\zL$ is a homogeneous (linear)
2-contravariant tensor field on $E$, i.e. $\zL$ is homogeneous of degree
-1 with respect to the Liouville vector field $\zD_E$, then
$<\zL,\D\zi_{\zm}\ot\D\zi_{\zn}>=\{\zi_\zm,\zi_\zn\}_\zL$ is again a linear
function associated with  an  element  $[\zm,\zn]_\zL$.  The  operation
$[\zm,\zn]_\zL$ on sections  of  $E^*$
we call the {\it bracket induced by $\zL$}.
This is the way in which homogeneous Poisson brackets are related
to Lie algebroids.

\begin{theorem} There is a one-one correspondence between  Lie  algebroid
brackets $\lna\cdot,\cdot\rna_\zL$ on the vector bundle $E$ and  homogeneous
(linear) Poisson structures $\zL$ on the dual bundle $E^*$ determined by
\be\label{lpo}
\zi_{\lna X,Y\rna_\zL}=\{\zi_X,\zi_Y\}_\zL=\zL(d\zi_X,d\zi_Y).
\ee
\end{theorem}

For  a  vector  bundle   $E$   over   the   base   manifold   $M$,   let
$\A(E)=\oplus_{k\in\Z}\A^k(E)$,   $\A^k(E)=Sec(\bigwedge^kE)$,    be    the
exterior
algebra of multisections of $E$.  This  is  a  basic  geometric model  for  a
graded
associative commutative algebra with unity. We will refer to elements of
$\zW^k(E)=\A^k(E^*)$ as to {\it k-forms} on $E$. Here, we identify
$\A^0(E)=\zW^0(E)$
with  the  algebra  $C^\infty(M)$  of smooth functions on the base and
$\A^k(E)=\{ 0\}$
for $k<0$. Denote by $\vert X\vert$ the Grassmann degree of the multisection
$X\in\A(E)$.

A Lie algebroid structure on $E$ can be
identified  with  a {\it graded Poisson bracket} on $\A(E)$ of
degree -1 (linear). Such brackets we call {\it Schouten-Nijenhuis brackets}
on $\A(E)$. Recall
that a graded Poisson bracket of degree  $k$  on  a  $\Z$-graded associative
commutative algebra $\A=\oplus_{i\in\Z}\A^i$ is a graded bilinear map
$$
\{\cdot,\cdot\}:\A\ti\A\ra\A
$$
of degree $k$ (i.e. $\vert\{a,b\}\vert=\vert
a\vert+\vert  b\vert+k$) such that
\begin{enumerate}
\item $\{ a,b\}=-(-1)^{(\vert a\vert+k)(\vert b\vert+k)}\{ b,a\}$ (graded
anticommutativity),
\item $\{ a,bc\}=\{   a,b\}   c+(-1)^{(\vert  a\vert+k)\vert  b\vert}b\{
a,c\}$  (graded Leibniz rule),
\item      $\{\{      a,b\},c\}=\{      a,\{        b,c\}\}-(-1)^{(\vert
a\vert+k)(\vert b\vert+k)}\{ b,\{ a,c\}\}$
(graded Jacobi identity).
\end{enumerate}

\noindent
It is obvious  that  this  notion  extends  naturally  to  more  general
gradings in the algebra. For  a  graded  commutative  algebra   with
unity   $\1$,   a natural
generalization of  a  graded  Poisson  bracket  is  {\it  graded  Jacobi
bracket}. The only difference is that we replace the Leibniz rule by the
{\it generalized Leibniz rule}
\be\label{first} \{ a,bc\}=\{
a,b\}  c+(-1)^{(\vert a\vert+k)\vert b\vert}b\{  a,c\}-\{ a,\1\} bc.
\ee
Graded Jacobi
brackets on $\A(E)$  of  degree  -1  (linear) we  call  {\it Schouten-Jacobi}
brackets.
An element $X\in \A^2(E)$  is  called  a  {\it canonical structure} for a
Schouten-Nijenhuis or Schouten-Jacobi bracket $\lna\cdot,\cdot\rna$ if $\lna
X,X\rna=0$.

As it was already indicated in \cite{KS},  Schouten-Nijenhuis brackets are in
one-one
correspondence with Lie algebroids:

\begin{theorem} Any  Schouten-Nijenhuis  bracket $\lna\cdot,\cdot\rna$  on
$\A(E)$ induces a Lie  algebroid  bracket   on   $\A^1(E)=Sec(E)$   with
the anchor
defined by $\zr(X)(f)=\lna X,f\rna$. Conversely, any Lie algebroid structure
on
$Sec(E)$ gives rise to a  Schouten-Nijenhuis  bracket  on $\A(E)$  for  which
$\A^1(E)=Sec(E)$  is  a  Lie  subalgebra  and $\zr(X)(f)=\lna X,f\rna$.
\end{theorem}
We have the following expression for the Schouten-Nijenhuis bracket:
\begin{eqnarray}\label{Schouten3}
\lefteqn
{\lna X_1\wedge\ldots\wedge X_m,Y_1\wedge\cdots\wedge Y_n\rna =} \\
& & \sum_{k,l}(-1)^{k+l}\lna X_k,Y_l\rna\wedge\ldots\wedge\widehat{X_k}\wedge
\ldots\we
X_m\we Y_1\we\ldots\we\widehat{Y_l}\we\ldots\we Y_n, \nonumber
\end{eqnarray}
where $X_i,Y_j\in Sec(E)$ and the hat over a symbol means that
this is to be omitted.

A Schouten-Nijenhuis bracket induces  the  well-known generalization of the
standard
Cartan calculus of differential forms and vector fields \cite{Ma,MX}. The
exterior
derivative $\D:\zW^k(E) \rightarrow\zW^{k+1}(E)$ is defined by the standard
formula
\bea\nn
        \D\zm(X_1,\dots,X_{k+1}) &=& \sum_i (-1)^{i+1}
\lna X_i,\zm(X_1,\dots,\widehat{X}_i,\dots ,X_{k+1})\rna \\
&+& \sum_{i<j} (-1)^{i+j}\zm (\lna X_i,X_j\rna, X_1,\dots , \widehat{X}_i,
\dots
,\widehat{X}_j, \dots , X_{k+1}),\label{D} \eea
        where $X_i\in Sec(E)$.  For $X \in Sec(E)$, the
contraction $i_X \colon \zW^p(E) \rightarrow\zW^{p-1}(E)$ is defined in the
standard
way and the Lie differential operator $ \Ll_X$
         is defined by the graded commutator
\be\label{L} \Ll_X = i_X \circ \D +\D \circ i_X. \ee

Since Schouten-Nijenhuis brackets  on  $\A(E)$  are  just  Lie
algebroid structures on $E$, by {\it Jacobi algebroid} structure
on $E$ we  mean  a Schouten-Jacobi bracket on $\A(E)$ (see
\cite{GM}). An  analogous concept has   been introduced in
\cite{IM} under the  name of  a  {\it  generalized  Lie
algebroid}. Every Schouten-Jacobi bracket on the graded algebra
$\A(E)$ of multisections of $E$ turns out to be uniquely
determined by a Lie algebroid bracket on a vector bundle $E$ over
$M$ and a 1-cocycle $\zF\in\zW^1(E)$, $\D\zF=0$, relative to the
Lie algebroid exterior derivative $\D$, namely it is of the form
\cite{IM} \be\label{jb} \lna X,Y\rnaf=\lna X,Y\rna+xX\we i_\zF
Y-(-1)^xyi_\zF X\we Y , \ee where $\lna,\cdot,\cdot\rna$ is the
Schouten bracket associated with this Lie algebroid and where we
use the convention  that  $x=\vert X\vert-1$ is  the  shifted
degree of $X$ in the graded algebra $\A(E)$. Note that $\zF$ is
determined by the Schouten-Jacobi bracket by $i_\zF X=(-1)^x\lna
X,\1\rnaf$, so that (\ref{first}) is satisfied: \be\label{first1}
\lna X,Y\we Z\rnaf=\lna X,Y\rnaf\we Z+(-1)^{x(y+1)}Y\we\lna
X,Z\rnaf - \lna X,\1\rnaf\we Y\we Z. \ee We already know that
there is one-one correspondence between Lie algebroid structures
on $E$  and linear Poisson  tensors  $\zL^{E^*}$  on  $E^*$.  To
Jacobi algebroids correspond  Jacobi  structures  $J_\zF^{E^*}$ on
$E^*$   which are homogeneous of degree -1 with respect  to the
Liouville  vector  field $\zD_{E^*}$, namely
$$  J_\zF^{E^*}=\zL^{E^*}+\zD_{E^*}\we\zF^v-I\we\zF^v,
$$
where $\zF^v$ is the vertical lift of $\zF$ to a vector field on $E^*$.
The above structure
generates  a   Jacobi   bracket   which   coincides on linear functions with
the
Poisson bracket associated with $\zL^{E^*}$.

One can develop a Cartan calculus for Jacobi algebroids similarly to the Lie
algebroid
case  (cf.  \cite{IM}).  For  a  Schouten-Jacobi  bracket associated with  a
1-cocycle  $\zF$  the  definitions  of  the  exterior differential  $\D^\zF$
and
Lie   differential   $\Ll^\zF=\D^\zF\circ i+i\circ\D^\zF$  are formally the
same as
(\ref{D}) and (\ref{L}),   respectively.     Since,     for   $X\in   Sec(E)$,
$f\in
C^\infty(M)$,   we have $\lna X,f\rnaf=\lna X,f\rna+(i_\zF X)f$, one obtains
$\D^\zF\zm=\D\zm+ \zF\we\zm$. Here  $\lna\cdot,\cdot\rna$  and  $\D$ are,
respectively,  the  Schouten-Nijenhuis  bracket  and  the   exterior
derivative
associated with the Lie algebroid.

\medskip\noindent
{\bf Example 1.} A canonical example of a Lie algebroid over $M$ is the
tangent
bundle
$TM$ with the bracket of vector fields. The corresponding complex
$(\zW(TM),\D)$ is in
this case the standard de Rham complex. A canonical structure  for the
corresponding
Schouten-Nijenhuis bracket is just a standard  Poisson tensor.

\medskip\noindent
{\bf Example 2.} A canonical example  of  a Jacobi algebroid is
$(T_1M=TM\oplus\R,(0,1))$,
where $T_1M$ is  the  Lie algebroid of first-order differential operators
on
$C^\infty(M)$  with the bracket
$$  [(X,f),(Y,g)]_1=([X,Y],X(g)-Y(f)),\quad
X,Y\in
Sec(TM),\quad   f,g\in C^\infty(M),
$$
and  the  1-cocycle   $\zF=(0,1)$   is
$\zF((X,f))=f$.   A   canonical structure with respect to the corresponding
Schouten-Jacobi  bracket  on the Grassmann algebra  $\A(T_1M)$  of
first-order
polydifferential operators on $C^\infty(M)$ turns out to be a standard Jacobi
structure. Indeed, it is  easy  to  see that the Schouten-Jacobi bracket reads
\bea\nn
&&\lna A_1+I\we  A_2,B_1+I\we  B_2\rna_1=\lna  A_1,B_1\rna  +(-1)^aI\we
\lna A_1,B_2\rna+I\we\lna A_2,B_1\rna\\
&&\qquad+aA_1\we B_2-(-1)^abA_2\we B_1+(a-b)I\we A_2\we B_2.\label{dif} \eea
Hence,
the bracket  $\{\cdot,\cdot\}$  on  $C^\infty(M)$  defined  by  a bilinear
differential operator $\zL+I\we\zG\in\A(T_1M)$  is  a  Lie bracket
(Jacobi
bracket on $C^\infty(M)$) if and only if
$$
\lna\zL+I\we\zG,\zL+I\we\zG\rna=\lna\zL,\zL\rna+2I\we
\lna\zG,\zL\rna+2\zL\we\zG=0.
$$
We recognize the conditions (\ref{kk}) defining a Jacobi structure on $M$.

There is another  approach  to  Lie  algebroids.  As  it has been
shown in \cite{GU1,GU2}, a Lie algebroid   structure   (or   the
corresponding Schouten-Nijenhuis bracket)  is  determined  by  the
Lie algebroid   lift $X\mapsto X^c$ which associates  with
$X\in\zT(E)$   a  contravariant tensor  field $X^c$ on $E$. The
complete Lie algebroid and  Jacobi  algebroid lifts     are
described  as follows.

\begin{theorem} (\cite{GU1}) For a given Lie  algebroid  structure  on  a
vector bundle $E$ over $M$ there is a unique {\it complete lift} of elements
$X\in Sec(E^{\ot k})$ of
the tensor algebra $\zT(E)=\oplus_k  Sec(E^{\ot k})$  to linear contravariant
tensors
$X^c\in Sec((TE)^{\ot k})$ on $E$ , such that
\begin{description}
\item{(a)} $f^c=\zi_{\D f}$ for $f\in C^\infty(M)$;
\item{(b)}  $X^c(\zi_\zm)=\zi_{\Ll_X\zm}$  for  $X\in Sec(E),\  \zm\in
Sec(E^*)$;
\item{(c)} $(X\ot Y)^c=X^c\ot Y^v+X^v\ot Y^c$, where $X\mapsto  X^v$  is
the standard vertical lift of tensors  from  $\zT(E)$  to  tensors  from
$\zT(TE)$, i.e. the complete lift is a derivation with  respect  to  the
vertical lift.
\end{description}
This  complete  lift restricted to skew-symmetric tensors  is a
homomorphism of   the corresponding Schouten-Nijenhuis brackets:
\be\label{l1}
\lna
X,Y\rna^c=\lna X^c,Y^c\rna.
\ee
Moreover,
\be\label{l2} \lna X,Y\rna^v=\lna X^c,Y^v\rna.
\ee
\end{theorem}
\begin{corollary} If $P\in\A^2(E)$ is  a  canonical  structure  for  the
Schouten bracket, i.e.   $\lna   P,P\rna=0$,  then  $P^c$
is a homogeneous Poisson structure on $E$.  The corresponding Poisson
bracket determines  the Lie algebroid bracket
\be\label{Fu0]}
\lna\za,\zb\rna_P=(i_{\#_P(\alpha)}\D\beta-
i_{\#_P(\beta)}\D\alpha +\D(P(\alpha,\beta))
\ee
on $E^*$.
\end{corollary}
\noindent {\bf Remark.} For the canonical Lie algebroid $E=TM$,
the above complete lift gives the better-known {\it tangent lift}
of multivector  fields on $M$ to multivector fields on $TM$ (cf.
\cite{IY,GU}). In this case the complete lift is an injective
operator, so $\zL$ is a Poisson tensor on $M$ if and only if
$\zL^c$ is a Poisson tensor on $TM$. The  complete Lie algebroid
lift of just sections of $E$, i.e. the  formula  (b),  was already
indicated in \cite{MX1}.

\medskip
Let us see how these lifts look like in local coordinates. Let
$(x^a)$ be a local coordinate system on $M$ and let $e_1,\dots,
e_n$ be  a basis of local sections of $E$. We denote by $
e^{*1},\dots, e^{*n}$ the dual basis of local sections of $E^*$
and by $(x^a,y^i)$ (resp. $(x^a,\zx_i)$) the corresponding
coordinate system on $E$ (resp. $E^*$), i.e., $\zi_{e_i} =\zx_i$
and $\zi_{e^{*i}}=y^i$. The vertical lift is given by
$$  (c_{i_1,\dots,i_k}e_{i_1}\ot\cdots\ot
e_{i_k})^v=
c_{i_1,\dots,i_k}\pa_{y^{i_1}}\ot\cdots\ot\pa_{y^{i_k}}.
$$
    If  for the Lie algebroid bracket we have $[e_i,e_j]=c_{ij}^ke_k$ and if
the anchor
sends $e_i$ to $d^a_i\pa_{x^a}$, then
\begin{equation}\label{Fp10}
\zL^{E^*} = \frac{1}{2}c^k_{ij}\zx_k \partial _{\zx_i}\wedge
\partial _{\zx_j} + d^a_i\partial _{\zx_i} \wedge \partial _{x^a}.
\end{equation}
Moreover,
\begin{equation}\label{Fp11}
f^c = \frac{\partial f}{\partial x^a} d^a_j y^{j}
\end{equation}
and
\begin{equation}\label{Fp12}
(X^ie_{i})^c = X^i d^a_i\partial _{x^a} + (X^i c^k_{ji} + \frac{\partial
X^k}{\partial
x^a} d^a_j )y^j \partial _{y^k}.
\end{equation}
It follows that, for $P=\frac{1}{2}P^{ij} e_i\we e_j $, we have
\begin{equation}\label{Fp13}
P^c = P^{ij}d^a_j \partial _{y^i}\we\partial _{x^a} + (P^{kj} c^i_{lk} +
\frac{1}{2}\frac{\partial P^{ij}}{\partial x^a} d^a_l)y^l \partial _{y^i}
\wedge
\partial _{y^j}.
\end{equation}

\medskip\noindent
There is an analog of the  Lie  algebroid  complete  lift  for
Jacobi algebroids which will represent the  Schouten-Jacobi
bracket  on $\A(E)$    in    the Schouten-Jacobi bracket     of
first-order polydifferential operators on $E$. Here by   {\it
polydifferential   operators}   we    understand    skew-symmetric
multidifferential operators. Let $\lna\cdot,\cdot\rnaf$ be the
Schouten-Jacobi bracket on  $\A(E)$ associated with   a  Lie
algebroid  structure  on  $E$  and  a 1-cocycle  $\zF$.

\medskip\noindent
{\bf Definition.}(\cite{GM}) The  {\it complete Jacobi  lift}  of
an element
$X\in\zT^k(E)$ is the multidifferential operator of first order on $E$,
i.e. an element of $Sec((T_1E)^{\ot k})$, defined by \be\label{jl} \widehat
X_\zF=X^c-(k-1)\zi_\zF X^v+i_{I\ot\D(\zi_\zF)} X^v, \ee where $X^c$ is the
complete
Lie algebroid lift, $X^v$ is the  vertical lift   and   $i_{I\ot\D(\zi_\zF)}$
is the
derivation    acting    on the tensor algebra  of  contravariant  tensor
fields
which  vanishes  on functions and satisfies $i_{I\ot\D(\zi_\zF)}X=X(\zi_\zF)I$
on vector fields. The derivation property yields
$$ i_{I\ot\D(\zi_\zF)}(X_1^v\ot\cdots\ot X_k^v)=
\sum_i\langle X_i,\zF\rangle X_1^v\ot\cdots\ot X^v_{i-1}\ot
I\ot X^v_{i+1}\ot\cdots\ot X^v_k.
$$
for $X_1\dots,X_k\in Sec(E)$.
\begin{theorem}(\cite{GM}) The complete Jacobi lift has the following
properties:
\begin{description}
\item{(a)} $\widehat f_\zF=\zi_{\D^\zF f}$ for $f\in C^\infty(M)$;
\item{(b)}   $\widehat X_\zF=X^c+(i_\zF X)^vI$
for $X\in Sec(E)$;
\item{(c)} $\widehat{(X\ot Y)}_\zF=\widehat X_\zF\ot Y^v+X^v\ot
\widehat Y_\zF-\zi_\zF(X^v\ot Y^v)$;
\item{(d)} For skew-symmetric tensors $X$ and $Y$,
$$  \lna\widehat X_\zF,\widehat Y_\zF\rna_1={(\lna X,Y\rnaf)}^\wedge_\zF,
$$
where $\lna\cdot,\cdot\rna_1$  is   the   Schouten-Jacobi
bracket   of first-order polydifferential operators;
\item{(e)} For skew-symmetric $X$ and $Y$
$$  \lna\widehat X_\zF,Y^v\rna_1=(\lna X,Y\rnaf)^v;
$$
\end{description}
\end{theorem}
Remark  that  in  \cite{GM}  only  skew-symmetric  tensors   have   been
considered, but the extension to arbitrary tensors is straightforward.

\medskip\noindent
{\bf Definition.}(\cite{GM}) The  {\it complete Poisson  lift}  of  an element
$X\in\zT^k(E)$  is  the   contravariant   tensor   $\widehat   X^c_\zF\in
Sec((TE)^{\ot k})$,
defined by \be\label{jl1} \widehat X^c_\zF=X^c-(k-1)\zi_\zF
X^v+i_{\zD_E\ot\D(\zi_\zF)}
X^v, \ee where $\zD_E$ is the Liouville vector field on the vector bundle $E$
and
$X^c$ is the complete Lie algebroid lift, $X^v$ is the vertical lift and
$i_{\zD_E\ot\D(\zi_\zF)}$ is the derivation acting on the  tensor algebra of
contravariant tensor fields which vanishes on  functions  and satisfies
$i_{\zD\ot\D(\zi_\zF)}X=X(\zi_\zF)\zD_E$ on vector fields.
\begin{theorem}(\cite{GM}) The Poisson lift has the following properties:
\begin{description}
\item{(a)} $\widehat f^c_\zF=\zi_{\D^\zF f}$ for $f\in C^\infty(M)$;
\item{(b)}   $\widehat X^c_\zF(\zi_\zm)=\zi_{\Ll^\zF_X\zm}$
for $X\in Sec(E)$, $\zm\in Sec(E^*)$;
\item{(c)} $\widehat{(X\ot Y)}^c_\zF=\widehat X^c_\zF\ot Y^v+X^v\ot
\widehat Y^c_\zF-\zi_\zF(X^v\ot Y^v)$;
\item{(d)} For skew-symmetric $X$ and $Y$
$$  \lna\widehat X^c_\zF,\widehat Y^c_\zF\rna={(\lna
X,Y\rnaf)}^{\wedge c}_\zF.
$$
\end{description}
\end{theorem}
\begin{corollary} If $P\in\A^2(E)$ is  a  canonical  structure  for  the
Schouten-Jacobi  bracket,   i.e.   $\lna   P,P\rnaf=\lna
P,P\rna+2P\we i_{\zF} P=0$,  then  $\widehat  P_{\zF}$ (resp.
$\widehat  P^c_{\zF}$)  is  a homogeneous Jacobi (resp.
homogeneous Poisson) structure on $E$.  The corresponding  Jacobi
and Poisson brackets coincide on linear functions and determine
the Lie algebroid bracket \be\label{Fu1}
\lna\za,\zb\rna_P=(i_{\#_P(\alpha)}\Df\beta-
i_{\#_P(\beta)}\Df\alpha +\Df(P(\alpha,\beta)) \ee on $E^*$.
\end{corollary}

\section{Characterization of Poisson tensors.}

 Theorem~1 of Introduction can be generalized in the  following way.  Let  us
remark first that any two-contravariant tensor   $\zL$ (which is not
assumed to be skew-symmetric)
defines a bracket $[\cdot,\cdot]_\zL$ on 1-forms on $M$ by \be\label{1}
[\zm,\zn]_\zL=i_{\sharp_\zL(\zm)}\D\zn-i_{\sharp_\zL(\zn)}\D\zm+
\D<\zL,\zm\otimes\zn>, \ee where $<\cdot,\cdot>$ is the canonical pairing
between
contravariant and covariant tensors.

\begin{theorem} For a two-contravariant tensor $\zL$ on a  manifold  $M$
the following are equivalent:
\begin{description}
\item{(i)} $\zL$ is a Poisson tensor;
\item{(ii)} $\sharp_\zL$ induces a homomorphism  of  $[\cdot,\cdot]_\zL$
into the bracket of vector fields:
\be\label{2}
\sharp_\zL([\zm,\zn]_\zL)=[\sharp_\zL(\zm),\sharp_\zL(\zn)];
\ee
\item{(iii)} The canonical Poisson tensor $\zL_M$ and the negative of the
complete lift $-\zL^c$ are $\sharp_\zL$-related;
\item{(iv)} There is a vector bundle morphism $F:T^*M\ra TM$  over  the
identity
on $M$ such that the canonical Poisson tensor $\zL_M$
and the negative of the complete lift $-\zL^c$ are $F$-related;
\item{(v)} The morphism $\sharp_\zL$ relates $\zL_M$ with  the  complete
lift of a 2-contravariant tensor $\zL_1$.
\item{(vi)} There is a vector bundle morphism $F:T^*M\ra TM$  over  the
identity
on $M$ such that
\be \label{2b}F([\zm,\zn]_\zL)=[F(\zm),F(\zn)].
\ee
\item{(vii)} There is a 2-contravariant tensor $\zL_1$ on $M$ such that
\be\label{2c}
\sharp_\zL([\zm,\zn]_{\zL_1})=[\sharp_\zL(\zm),\sharp_\zL(\zn)];
\ee
\end{description}\label{t2}
\end{theorem}

\begin{pf} The implication $(i)\Rightarrow(ii)$ is a well-known fact (cf.
e.g. \cite{KSM}).

Assume now $(ii)$. To show $(iii)$ one has to prove
that   the   brackets   on   functions   $\{\cdot,\cdot\}_{\zL_M}$    and
$\{\cdot,\cdot\}_{\zL^c}$ induced  by  tensors  $\zL_M$  and  $\zL^c$  by
contractions with differentials of functions are
$\sharp_\zL$-related, i.e.
    \begin{equation}\label{3}
-\{ f,g\}_{\zL^c}\circ\sharp_\zL=\{  f\circ\sharp_\zL,g\circ\sharp_\zL\}_
{\zL_M}
    \ee
    for all $f,g\in C^\infty(TM)$. Due to Leibniz rule, it is  sufficient
to check (\ref{3}) for linear functions, i.e. for functions of the
form $\zi_\zm$,  where $\zm$
is  a 1-form and  $\zi_\zm(v_x)=<\zm(x),v_x>$.
It   is   well   known   (see \cite{Co,GU}) that the brackets induced by
$\zL$  and  its  complete  lift are related by
\be\label{l} \{\zi_\zm,\zi_\zn\}_{\zL^c}=\zi_{[\zm,\zn]_\zL}.
\ee
It is also known (cf. (\ref{lpo}))
that
$$ \zi_{[X,Y]}=\{\zi_X,\zi_Y\}_{\zL_M}
$$
for vector fields $X,Y$ on $M$.
Since $\zi_\zm\circ\sharp_\zL=-\zi_{\sharp_\zL(\zm)}$, we get \bea\label{4}
-\{\zi_\zm,\zi_\zn\}_{\zL^c}\circ\sharp_\zL&=&-\zi_{[\zm,\zn]_\zL}
\circ\sharp_\zL=\zi_{\sharp_\zL([\zm,\zn]_\zL)}\\
&=&\zi_{[\sharp_\zL(\zm),\sharp_\zL(\zn)]}=
\{\zi_{\sharp_\zL(\zm)},\zi_{\sharp_\zL(\zn)}\}_{\zL_M}\nn \\
&=&\{\zi_{\zm}\circ\sharp_\zL,\zi_{\zn}\circ\sharp_\zL\}_{\zL_M}\nn
\eea
which proves $(ii)\Rightarrow(iii)$. In fact, (\ref{4}) proves equivalence
of $(ii)$ and $(iii)$.

Replacing in (\ref{4}) the mapping $\sharp_\zL$ by a vector bundle morphism
$F:T^*M\ra TM$, we get equivalence of $(iv)$ and $(vi)$.
Similarly, $(v)$ is equivalent to $(vii)$.
The implication $(iii)\Rightarrow(iv)$ is obvious, so let us show
$(iv)\Rightarrow(i)$. Assume that $F$ relates $\zL_M$  and  $\zL^c$.  We
will  show  that  this  implies  that  $\zL$   is   skew-symmetric   and
$F=\sharp_\zL$. Since the  assertion  is  local  over  $M$  we  can  use
coordinates  $(x^a)$  in  $M$  and  the   adapted   coordinate   systems
$(x^a,p_i)$  in  $T^*M$   and   $(x^a,\dot x^j)$   in   $TM$.   Writing
$\zL=\zL^{ij}\pa_{x^i}\otimes\pa_{x^j}$                              and
$F(x^a,p_i)=(x^a,F^{ij}p_i)$, we get
\be\label{m}
F_*(\pa_{p_i}\we\pa_{x^i})=F^{ij}p_s\frac{\pa F^{sk}}{\pa x^i}
\pa_{\dot x^j} \we\pa_{\dot x^k}-F^{ij}\pa_{x^i}\we\pa_{\dot x^j}.
\ee
Since
\be\label{lift}
\zL^c=\frac{\pa\zL^{ij}}{\pa x^k}\dot x^k\pa_{\dot x^i}\otimes \pa_{\dot
x^j}+\zL^{ij}(\pa_{x^i}\otimes\pa_{\dot  x^j}+\pa_{\dot     x^i}\otimes
\pa_{x^j}),
\ee
comparing the vertical-horizontal parts we get  $\zL^{ij}=F^{ij}=-F^{ji}$,
i.e. $\zL$ is skew-symmetric and $F=\sharp_\zL$.  Going  backwards  with
(\ref{4}) we get $(ii)$. But for  skew  tensors  we
have (cf. \cite{KSM})
\be\label{5}
\frac{1}{2}\lna\zL,\zL\rna(\zm,\zn,\zg)=<\sharp_\zL([\zm,\zn]_\zL)-
[\sharp_\zL(\zm),\sharp_\zL(\zn)],\zg>,
\ee
where $\lna\cdot,\cdot\rna$ is the Schouten-Nijenhuis  bracket,  so  that
$\lna\zL,\zL\rna=0$, i.e. $\zL$ is a Poisson tensor.

Finally, $(v)$   is equivalent  to  $(iii)$, since
$(iii)\Rightarrow(v)$  trivially and exchanging the role of $F$
and $\zL$ in (\ref{m}) and (\ref{lift}) we see that, as above,
$\zL^{ij}=F^{ij}$, so that any tensor whose complete lift is
$\sharp_\zL$-related to $\zL_M$ equals $-\zL$.
\end{pf}

\medskip
A  similar  characterization  is  valid  for  any  Lie algebroid.
Let us consider a vector bundle $E$ over $M$ with a Lie algebroid
bracket $\lna\cdot,\cdot\rna$ instead of the canonical Lie
algebroid $TM$  of  vector fields (cf. \cite{Ma,KSM,GU0}). The
multivector fields are now  replaced by multisections
$\A(E)=\oplus_k\A^k(E)$, $\A^k(E)=Sec(\bigwedge^kE)$, of $E$ and
the standard Schouten-Nijenhuis bracket  with  its Lie algebroid
counterpart. A Lie algebroid Poisson tensor ({\it canonical
structure}) is then a skew-symmetric  $\zL\in \A^2(E)$  satisfying
$\lna\zL,\zL\rna=0$. Such a structure gives a {\it triangular Lie
bialgebroid} in the sense  of \cite{MX}. We have the exterior
derivative $\D$  on multisections  of the dual bundle $E^*$ (we
will refer to them as to  "exterior  forms"). For any $\zL\in
Sec(E\otimes E)$  the formula  (\ref{1})  defines  a bracket on
"1-forms". We have  an analog  of  the  complete  lift  (cf.
\cite{GU1,GU2})
$$ Sec(E^{\otimes k})\ni \zL\mapsto\zL^c\in Sec((TE)^{\otimes k})
$$
of the tensor algebra of sections of $E$ into contravariant  tensors  on
the total space $E$.  The  Lie algebroid
bracket  corresponds  to  a  linear  Poisson   tensor   $\zL^{E^*}$   on
$E^*$ (which is just $\zL_M$ in the case $E=TM$) by (\ref{lpo}).
Since the tensor  $\zL^{E^*}$  may  be  strongly  degenerate,
linear maps $F:E^*\ra E$ do not determine the related tensors  uniquely,
so we cannot have the full analog of Theorem \ref{t2}. However, since for
skew-symmetric tensors the formula (\ref{5}) remains valid \cite{KSM}, a
part of Theorem \ref{l2} can be proved in the same way, {\it mutatis
mutandis},
in the general Lie algebroid case. Thus we get the Lie algebroid version
of Theorem 1 (cf. \cite{GU1}).

\begin{theorem} For any bisection $\zL\in \A^2(E)$ of a Lie algebroid $E$
the following are equivalent:
\begin{description}
\item{(i)} $\zL$ is a canonical structure, i.e. $\lna \zL,\zL\rna=0$;
\item{(ii)} $\sharp_\zL$ induces a homomorphism  of  $[\cdot,\cdot]_\zL$
into the Lie algebroid bracket:
\be\label{2a}
\sharp_\zL([\zm,\zn]_\zL)=[\sharp_\zL(\zm),\sharp_\zL(\zn)];
\ee
\item{(iii)} The canonical Poisson tensor $\zL_M$ and the negative of the
complete lift $-\zL^c$ are $\sharp_\zL$-related.
\end{description}\label{t3}
\end{theorem}
\section{Jacobi algebroids and characterization of Jacobi structures}

We have introduced in Section~2 Jacobi and Poisson complete lifts related to
Jacobi algebroids.
For  a standard Jacobi structure $J=(\zL,\zG)$  on $M$  we  will  denote
these lifts of $J$ by $\widehat J$ and $\widehat J^c$, respectively. The
Jacobi lift
$\widehat J$ is the Jacobi structure on $E=TM\oplus\R$ given by \cite{GM}
\be\label{im}
\widehat J=(\zL^c-t\zL^v+\pa_t\we(\zG^c-t\zG^v),\zG^v),
\ee
where $\zL^v$ and $\zG^v$ are the vertical tangent lifts  of  $\zL$ and
$\zG$, respectively, and $t$ is the standard linear coordinate in $\R$.
We consider here tangent lifts as tensors on
$TM\oplus\R=TM\ti\R$ instead on $TM$. The linear Jacobi structure (\ref{im})
has been  already  considered by Iglesias and Marrero \cite{IM1}.

Similarly, the Poisson  lift  $\widehat J^c$  is  the  linear  Poisson
tensor  on $TM\oplus\R$ given by \cite{GM}
\be\label{gm}
\widehat J^c=\zL^c-t\zL^v+\pa_t\we\zG^c+\zD_{TM}\we\zG^v,
\ee
where $\zD_{TM}$ is the Liouville (Euler) vector  field  on  the  vector
bundle $TM$. This is exactly the linear Poisson tensor corresponding  to
the  Lie  algebroid  structure  on  $T^*M\oplus\R$  induced  by  $J$  and
discovered first in \cite{KSB}:
\bea\nn
[(\za,f),(\zb,g)]_J&=&(\Li_{\sharp_\zL(\za)}\zb-\Li_{\sharp_\zL(\zb)}\za
-\D<\zL,\za\we\zb>
+f\Li_\zG\zb-g\Li_\zG\za-i_\zG\za\we\zb,\label{ks}\\ \label{ks1}
&&<\zL,\zb\we\za>+\sharp_\zL(\za)(g)-\sharp_\zL(\zb)(f)+f\zG(g)-g\zG(f)),
\eea
Of course, these lifts and an analog of the bracket (\ref{ks1}) are
well-defined for any first-order bidifferential operator
\be\label{do}
J=\zL+I\otimes\zG_1+\zG_2\otimes I+\za I\otimes I,
\ee
where $\zL$ is a 2-contravariant tensor, $\zG_1,\zG_2$ are vector fields,
and $\za$ is a function on $M$. The associated bracket
acts  on functions on $M$ by
$$ \{ f,g\}_J=<\zL,\D f\otimes\D g>+f\zG_1(g)+g\zG_2(f)+\za fg,
$$
The Jacobi  lift  of  $J$  is  the  first-order
bidifferential operator on $TM\oplus\R$ given by
\beas
\widehat J&=&\zL^c-t\zL^v+\pa_t\otimes(\zG_1^c-t\zG_1^v)+(\zG_2^c-t\zG_2^v)
\otimes\pa_t+(\za^c-t\za^v)\pa_t\otimes\pa_t\\
&&+I\otimes(\zG_1^v+\za^v\pa_t)+(\zG_2^v+\za^v\pa_t)\otimes I
\eeas
and the Poisson lift is the 2-contravariant tensor field
\beas
\widehat J^c&=&\zL^c-t\zL^v+\pa_t\otimes(\zG_1^c -t\zG_1^v)+(\zG_2^c -
t\zG_2^v)
\otimes\pa_t+(\za^c-t\za^v)\pa_t\otimes\pa_t\\
&&+\zD_{E}\otimes(\zG_1^v+\za^v\pa_t)+
(\zG_2^v+\za^v\pa_t)\otimes\zD_{E}\\
&=&\zL^c-t\zL^v+\pa_t\otimes\zG_1^c +\zG_2^c
\otimes\pa_t+(\za^c+t\za^v)\pa_t\otimes\pa_t\\
&&+\zD_{TM}\otimes\zG_1^v+ \zG_2^v\otimes\zD_{TM}. \eeas The
mapping $\sharp_J:E^*=T^*M\oplus\R\ra E=TM\oplus\R$ reads
$$
\sharp_J(\zw_x,\zl)=(\sharp_\zL(\zw_x)+\zl\zG_1(x),\zG_2(x)(\zw_x)+\za(x
)\zl).
$$
Note that any morphism from the vector bundle $E^*=T^*M\oplus\R$
into  $E=TM\oplus\R$ over the identity on $M$ is of this form.

The bidifferential operators $\widehat J$ and $\widehat J^c$ define brackets
$\{\cdot,\cdot\}_{\widehat J}$ and $\{\cdot,\cdot\}_{\widehat J^c}$,
respectively, on
functions on $TM\oplus\R$. These brackets coincide on linear
functions which close on a subalgebra with respect to them, so that they
define the bracket $[\cdot,\cdot]_J$ on sections of $T^*M\oplus\R$
(which coincides with the bracket (\ref{ks}) for skew-symmetric operators) by
$$ \{\zi_{(\zm,f)},\zi_{(\zn,g)}\}_{\widehat J}=
\{\zi_{(\zm,f)},\zi_{(\zn,g)}\}_{J^c}=\zi_{[(\zm,f),(\zn,g)]_J},
$$
where  $\zm,\zn$  are  1-forms,  $f,g$  are  functions   on   $M$,   and
$\zi_{(\zm,f)}=\zi_\zm+tf^v$.
Here we identify $T^*M\oplus\R$ with  $T^*M\ti\R$  and  use  the  linear
coordinate $\zl$ in $\R$. For the similar identification of $TM\oplus\R$
we use the coordinate $t$ of $\R$ in $TM\ti\R$, since  both  $\R$'s  play
dual roles.

  We  have  two canonical structures on the vector bundle
$E^*=T^*M\oplus\R\simeq  T^*M\ti\R$. One is the Jacobi structure (bracket)
\be\label{jm}
J_M=\zL_M+\zD_{T^*M}\we\pa_\zl+\pa_\zl\we I
\ee
and the other is the Poisson structure $\zL_M$ regarded as  the  product
of  $\zL_M$  on  $T^*M$  with  the  trivial  structure  on  $\R$.  These
brackets coincide  on  linear  functions which  close  on  a  subalgebra
with respect to both brackets,  so  that  they
define a Lie algebroid  structure  on the dual bundle  $E=TM\oplus\R$.
This  is the  Lie algebroid of first-order differential operators with
the bracket
$$ [(X,f),(Y,g)]_1=([X,Y],(X(g)-Y(f)),
$$
where $X,Y$ are vector fields and $f,g$ are functions on $M$.

\begin{theorem} For a first-order bidifferential operator $J$ the
following  are equivalent:
\begin{description}
\item{(J1)} $J$ is a Jacobi bracket;
\item{(J2)} The canonical Jacobi bracket $J_M$
and $-\widehat J$ are $\sharp_J$-related;
\item{(J3)}  There  is a first-order bidifferential operator $J_1$   such
that  $J_M$  and $-\widehat J$   are $\sharp_{J_1}$-related;
\item{(J4)}  There is a first-order bidifferential operator  $J_1$   such
that  $J_M$  and  $-\widehat  J_1$   are $\sharp_J$-related;
\item{(J5)}   The   contravariant   tensors   $\zL_M$   and   $-\widehat J^c$
are $\sharp_J$-related;
\item{(J6)}  There  is a first-order bidifferential operator $J_1$   such
that  $\zL_M$  and   $-\widehat J^c$   are $\sharp_{J_1}$-related;
\item{(J7)}  There  is a first-order bidifferential operator $J_1$   such
that  $\zL_M$  and  $-\widehat J_1^c$   are $\sharp_J$-related;
\item{(J8)} For any 1-forms $\zm,\zn$ and functions $f,g$ on $M$
$$ \sharp_J([(\zm,f),(\zn,g)]_J)=[\sharp_J(\zm,f),\sharp_J(\zn,g)]_1.
$$
\item{(J9)}  There  is a first-order bidifferential operator $J_1$   such
that
$$ \sharp_{J_1}([(\zm,f),(\zn,g)]_J)=[\sharp_{J_1}(\zm,f),
\sharp_{J_1}(\zn,g)]_1.
$$
\item{(J10)}  There  is a first-order bidifferential operator $J_1$   such
that
$$ \sharp_J([(\zm,f),(\zn,g)]_{J_1})=[\sharp_J(\zm,f),\sharp_J(\zn,g)]_1.
$$
\end{description}\label{t4}
\end{theorem}

Before proving this theorem we  introduce  some  notation  and  prove  a
lemma.
For a first-order bidifferential operator $J$ as in (\ref{do}), the {\it
poissonization} of $J$ is the tensor field on $M\ti\R$ of the form
\be\label{po}
P_J=e^{-s}(\zL+\pa_s\ot\zG_1+\zG_2\ot\pa_s+\za\pa_s\ot\pa_s),
\ee
where $s$ is the coordinate on $\R$.
Identifying $T^*(M\ti\R)$ with $T^*M\ti T^*\R$
(with coordinates $(s,\zl)$ in $T^*\R$) and $T(M\ti\R)$ with $TM\ti  T\R$
(with coordinates $(s,t)$ in $T\R$) we can write
\beas
\zL_{M\ti\R}&=&\zL_M+\pa_\zl\we\pa_s,\\
P_J^c&=&e^{-s}(\zL^c-t\zL^v+\pa_t\ot(\zG_1^c-t\zG_1^v)+
(\zG_2^c-t\zG_2^v)\ot\pa_t\\
&&+\pa_s\ot(\zG_1^v+\za^v\pa_t)+
(\zG_2^v+\za^v\pa_t)\ot\pa_s+(\za^c-t\za^v)\pa_t\ot\pa_t,\\
\sharp_{P_J}(\zw_x,\zl_s)&=&e^{-s}(\sharp_\zL(\zw_x)+\zl_s\zG_1(x),
\zG_2(x)(\zw_x)+\zl_s\za(x)).
\eeas
In local coordinates $x=(x^l)$ on  $M$  and  adapted  local  coordinates
$(x,p)$ on $T^*M$ and $(x,\dot x)$ on $TM$ we have
$$ (x^l,s,\dot
x^i,t)\circ\sharp_{P_J}=(x^l,s,e^{-s}(\zL^{ki}p_k+\zl\zG_1^i),
e^{-s}(\zG_2^kp_k+\zl\za))
$$
for    $\zL=\zL^{ij}\pa_{x^i}\ot\pa_{x^j}$,    $\zG_u=\zG_u^k\pa_{x^k}$,
$u=1,2$.
It is well known that $P_J$ is a Poisson  tensor  if
and only if $J$ is a Jacobi structure \cite{Li,GL}. In view of
Theorem \ref{t3}, we
can conclude  that  $J$  is  a  Jacobi  structure  if  and  only   if
$\zL_{M\ti\R}$  and  $P_J^c$  are  related  by  the  map  $\sharp_{P_J}:
T^*(M\ti\R)\ra T(M\ti\R)$.
Since $T(M\ti\R)\simeq E\ti\R$ and $T^*(M\ti\R)\simeq E^*\ti\R$, we  can
consider  the   bundles   $E=TM\oplus\R$   and   $E^*=T^*M\oplus\R$   as
submanifolds of $T(M\ti\R)$ and $T^*(M\ti\R)$,  respectively,  given  by
the equation $s=0$.

For any function $\zvf\in C^\infty(E)$ we  denote  by
$\breve\zvf$    the    function    on     $T(M\ti\R)=E\ti\R$     given     by
$\breve\zvf(v_x,s)=e^s\zvf(v_x)$.  Similarly,   for   any   function   $\zf\in
C^\infty(E^*)$ we denote by $\tilde\zf$ the function on
$T^*(M\ti\R)=E^*\ti\R$ given by
$\tilde\zf(u_x,s)=e^s\zf(e^{-s}u_x)$.

It is a matter of easy calculations to prove the following.
\begin{lemma}
\begin{description}
\item{(a)}       The        maps         $\zvf\mapsto \breve\zvf$         and
$\zvf\mapsto\tilde\zvf$  are
injective.
\item{(b)} For any first-order bidifferential operator $J$,
$$ \breve\zvf\circ\sharp_{P_J}={(\zvf\circ\sharp_J)}^\sim.
$$
\item{(c)} For any $\zvf,\zc\in C^\infty(E)$,
$$ \{\breve\zvf,\breve\zc\}_{P^c_J}=(\{\zvf,\zc\}_{\widehat J})^\smile.
$$
\item{(d)} For any $\zvf,\zc\in C^\infty(E^*)$,
$$ \{\tilde\zvf,\tilde\zc\}_{\zL_{M\ti\R}}=
{(\{\zvf,\zc\}_{J_M})}^\sim.
$$
\item{(e)} For linear $\zvf,\zc\in C^\infty(E)$,
$$ \{\breve\zvf,\breve\zc\}_{P^c_J}=(\{\zvf,\zc\}_{J^c})^\smile.
$$
\item{(f)} For linear $\zvf,\zc\in C^\infty(E^*)$,
$$ \{\tilde\zvf,\tilde\zc\}_{\zL_{M\ti\R}}=
{(\{\zvf,\zc\}_{\zL_M})}^\sim.
$$
\end{description}
\end{lemma}

\begin{pft}{\it\ref{t4}}. Due to the above Lemma  the following identities
are valid  for  arbitrary  $\zvf,  \zc\in   C^\infty(E)$   and
arbitrary first-order bidifferential operators $J,J_1$:
\beas
{(\{\zvf,\zc\}_{\widehat J}\circ\sharp_{J_1})}^\sim&=&
(\{\zvf,\zc\}_{\widehat J})^\smile\circ\sharp_{P_{J_1}}=
\{\breve\zvf,\breve\zc\}_{P^c_J}\circ\sharp_{P_{J_1}}\\
(\{\zvf\circ\sharp_{J_1},\zc\circ\sharp_{J_1}\}_{J_M})^\sim&=&
\{(\zvf\circ\sharp_{J_1})^\sim,(\zc\circ\sharp_{J_1})^\sim\}_
{\zL_{M
\ti\R}}=\{\breve\zvf\circ\sharp_{P_{J_1}},\breve\zc\circ\sharp_{P_{J_1}}
\}_{\zL_{M\ti\R}}.
\eeas
Thus
$$ -\{\zvf,\zc\}_{\widehat J}\circ\sharp_{J_1}=
\{\zvf\circ\sharp_{J_1},\zc\circ\sharp_{J_1}\}_{J_M}
$$
if and only if
\be\label{pss}
-\{\breve\zvf,\breve\zc\}_{P^c_J}\circ\sharp_{P_{J_1}}=
\{\breve\zvf\circ\sharp_{P_{J_1}},\breve\zc\circ\sharp_{P_{J_1}}
\}_{\zL_{M\ti\R}},
\ee
which  means  that $J_M$ and $-\widehat J$ are $\sharp_{J_1}$-related
if and only if  $\zL_{M\ti\R}$  and  the  complete   lift   of   the
poissonization $-P^c_J$ are $\sharp_{P_{J_1}}$-related. Due  to  Theorem
\ref{t3}, we get that $P_{J_1}=P_J$ and the poissonization $P_J$  is  a
Poisson tensor what, in turn, is equivalent to the fact that $J$  is  a
Jacobi bracket. Thus we get
$$
(J1)\Leftrightarrow(J2)\Leftrightarrow(J3)\Leftrightarrow(J4).$$
Using now linear functions $\zvf,\zc$, we get in  a
similar way that (\ref{pss}) is equivalent to
$$ -\{\zvf,\zc\}_{\widehat J^c}\circ\sharp_{J_1}=
\{\zvf\circ\sharp_{J_1},\zc\circ\sharp_{J_1}\}_{\zL_M}
$$
which, due to Theorem \ref{t3}, gives
$$
(J1)\Leftrightarrow(J5)\Leftrightarrow(J6)\Leftrightarrow(J7).$$
Finally,    completely    analogously    to     (\ref{4})     we     get
$(J5)\Leftrightarrow(J8)\Leftrightarrow(J9)\Leftrightarrow(J10)$.
\end{pft}

\medskip
{\bf Remark.} In the above proof we get the lifts $\widehat J$,  $J^c$,
and the map  $\sharp_J$  in a natural way  by  using the poissonization and
its
tangent lift.  This  is a geometric version of the methods in \cite{Va} for
obtaining $J^c$.
Note also that $J_M$  is  the  canonical  Jacobi  structure   on
$T^*M\ti\R$ regarded as a  contact  manifold in a natural way and  that  the
equivalence $(J1)\Leftrightarrow(J8)$ is a version of the characterization
in \cite{MMP}.

\medskip
The above theorem characterizing Jacobi structures one can generalize to
canonical structures associated with Jacobi algebroids as follows.

Consider now a Jacobi algebroid, i.e. a vector bundle  $E$  over
$M$ equipped with a  Lie  algebroid  bracket  $[\cdot,\cdot]$  and a `closed
1-form' $\zF\in\zW^1(E)$. We denote by  $\lna,\cdot,\cdot\rna$ the
Schouten-Nijenhuis bracket of the Lie algebroid and by  $\zT(E)\ni X\mapsto
X^c\in\zT(TE)$ the complete lift from the tensor algebra of $E$ into the
tensor algebra of $TE$. The corresponding Schouten-Jacobi
bracket  we  denote  by  $\lna\cdot,\cdot\rnaf$  and  the  corresponding
complete Jacobi and Poisson lifts by $\zT(E)\ni X\mapsto  \widehat  X_\zF\in
\zT(TE)$   and   $\zT(E)\ni   X\mapsto   X_\zF^c\in    \zT(TE)$,
respectively.

If the 1-cocycle $\zF$ is exact, $\zF=\D s$, we can obtain
the bracket $\lna\cdot,\cdot\rnaf$ from $\lna\cdot,\cdot\rna$ using  the
linear  automorphism  of   $\A(E)$   defined   by   $\A^k(E)\ni   X\mapsto
e^{-(k-1)s}X$ (cf. \cite{GM1}). This is a version of the Witten's  trick
\cite{Wi} to obtain the deformed  exterior  differential
$\D^\zF\zm=\D\zm+\zF\we\zm$ via  the
automorphism of the cotangent bundle given by multiplication by $e^s$.

Even if the 1-cocycle $\zF$ is not
exact, there is a nice construction \cite{IM} which allows to view $\zF$
as being exact but for an extended Lie algebroid in the bundle $\widehat
E=E\ti\R$ over $M\ti\R$. The sections of this bundle may be viewed as
parameter-dependent (s-dependent) sections of $E$. The sections of $E$ form a
Lie subalgebra
of s-independent  sections  in  the  Lie  algebroid  $\widehat  E$  which
generate the $C^\infty(M\ti\R)$-module of sections of $\widehat E$  and  the
whole structure is uniquely determined by putting the anchor $\widehat\zr(X)$
of a s-independent section $X$ to be $\widehat\zr(X)=\zr(X)+\langle\zF,X
\rangle\pa_s$,
where $s$ is the standard coordinate function in $\R$ and $\zr$ is
the anchor in $E$. All this is consistent (thanks to the fact that $\zF$
is a cocycle) and defines a Lie algebroid structure on $\widehat E$ with the
exterior derivative $\D$ satisfying $\D s=\zF$.

\medskip
Let now $U:\zT(E)\ra \zT(\widehat E)$ be natural embedding of the tensor
algebra of $E$ into the tensor subalgebra of s-independent sections
of $\widehat E$. It is obvious that on skew-symmetric tensors $U$ is a
homomorphism of the corresponding Schouten
brackets:
$$ \lna U(X),U(Y)\rnah=U(\lna X,Y\rna),
$$
where  we  use  the  notation  $\lna\cdot,\cdot\rna$   and   $\lna\cdot,
\cdot\rnah$ for the Schouten brackets in $E$ and $\widehat E$, respectively.
Let us now gauge $\zT(E)$ inside $\zT(\widehat E)$ by putting
$$ P^\zF(X)=e^{-ks}U(X)
$$
for any element $X\in Sec(E^{\ot (k+1)})$. Note that $X\mapsto
P^\zF(X)$Here by   {\it   polydifferential   operators}   we
understand    skew-symmetric multidifferential operators.
preserves  the grading but not the tensor product. It can be
easily proved (cf. \cite{GM1}) that the Schouten-Jacobi bracket
$\lna\cdot,\cdot\rnaf$ can  be obtained by this gauging from the
Lie algebroid bracket.
\begin{theorem}(\cite{GM1})  For  any  $X\in \A(E),
Y\in \A^(E)$ we have \be\label{j} \lna
P^\zF(X),P^\zF(Y)\rnah=P^\zF({\lna X,Y\rnaf}). \ee
\end{theorem}
We will usually skip the symbol $U$ and write simply $P^\zF(X)=e^{-ks}X$,
regarding $\zT(E)$ as embedded in $\zT(\widehat E)$. The complete  lift  for
the Lie algebroid $\widehat E$ will be denoted by $X\mapsto X^{\widehat  c}$
to distinguish from the lift for $E$. It is easy to see that
$$ (P^\zF(X))^{\widehat c}=(e^{-ks}X)^{\widehat  c}=e^{-ks}(X^c-k\zi_\zF
X^v+\pa_s\we
(i_\zF X)^v).
$$
Here we understand tensors on $E$ as tensors on $\widehat  E=E\ti\R$  in
obvious way.
Note  that  $\widehat{(E^*)}=(\widehat
E)^*$ and the linear Poisson tensor $\zL^{\widehat E^*}$ reads
$$ \zL^{\widehat E^*}=\zL^{E^*}+\zF^v\we\pa_s,
$$
where $\zL^{E^*}$  is  the  Poisson  tensor  corresponding  to  the  Lie
algebroid $E$ and $\zF^v$ is the vertical lift of $\zF$. Recall that  on
$E^*$ we  have also  a canonical Jacobi structure
$$ J_\zF^{E^*}=\zL^{E^*}+\zD_{E^*}\we\zF^v-I\we\zF^v
$$
which generates a  Jacobi  bracket  which  coincides  with  the  Poisson
bracket of $\zL^{E^*}$ on linear functions.

Let us remark that the map $P^\zF$  plays  the  role  of  a  generalized
poissonization.  Indeed,  for  the  Jacobi  algebroid   of   first-order
differential operators  $E=TM\oplus\R$ the extended Lie algebroid $\widehat
E\ti\R$ is  canonically isomorphic with $T(M\ti\R)$, $U((X,f))=X+f\pa_s$,
and for $J\in Sec(E^{\ot 2})$ the tensor field $P^\zF(J)$ coincides with
(\ref{po}).

Let now $J\in \A^2(E)$. The tensor $J$ is  a  canonical  structure   for
the Jacobi algebroid $(E,\zF)$,  i.e.  $\lna  J,J\rnaf=0$,  if  and  only  if
$P^\zF(J)$ is a
canonical structure for the Lie algebroid $\widehat E$, i.e. $\lna
P^\zF(J),P^\zF(J)\rnah=0$. Moreover,
$$
\sharp_{P^\zF(J)}(u_x,s)=(e^{-s}\sharp_J(u_x),s).
$$
Like above, for any function
$\zvf\in C^\infty(E)$ we  denote  by $\breve\zvf$    the    function    on
$\widehat
E=E\ti\R$     given     by $\zvf_E(v_x,s)=e^s\zvf(v_x)$ and   for   any
function
$\zf\in C^\infty(E^*)$ we denote by $\tilde\zf$ the function on $\widehat
E^*=E^*\ti\R$ given by $\tilde\zf(u_x,s)=e^s\zf(e^{-s}u_x)$. Recall that
(cf. Section 2)
$$  \widehat J_\zF=J^c-\zi_\zF  J^v+I\we(i_\zF  J)^v
$$
and
$$ \widehat
J_\zF^c=J^c-\zi_\zF J^v+\zD_E\we(i_\zF J)^v.
$$
The corresponding  brackets on
functions  on  $E$ coincide  on  linear functions and define a bracket
$[\cdot,\cdot]_J$ on sections of $E^*$  in the standard way:
$$
\zi_{[\zm,\zn]_J}=\{\zi_\zm,\zi_\zn\}_{\widehat J_\zF}
=\{\zi_\zm,\zi_\zn\}_{\widehat
J_\zF^c}
$$
Completely analogously to Lemma 1 we get the following.

\begin{lemma}
\begin{description}
\item{(a)}       The        maps         $\zvf\mapsto\breve\zvf$         and
$\zf\mapsto\tilde\zf$  are
injective.
\item{(b)} For any $J\in \A^2(E)$
$$ \breve\zvf\circ\sharp_{P^\zF(J)}={(\zvf\circ\sharp_J)}^\sim.
$$
\item{(c)} For any $\zvf,\zc\in C^\infty(E)$
$$ \{\breve\zvf,\breve\zc\}_{(P^\zF(J))^{\widehat
c}}=(\{\zvf,\zc\}_{\widehat J_\zF})^\smile.
$$
\item{(d)} For any $\zvf,\zc\in C^\infty(E^*)$
$$ \{\tilde\zvf,\tilde\zc\}_{\zL^{\widehat E^*}}=
{(\{\zvf,\zc\}_{J_\zF^{E^*}})}^\sim.
$$
\item{(e)} For linear $\zvf,\zc\in C^\infty(E)$
$$ \{\breve\zvf,                     \breve\zc\}_{(P^\zF(J))^{\widehat
c}}=(\{\zvf,\zc\}_{\widehat J^c})^\smile.
$$
\item{(f)} For linear $\zvf,\zc\in C^\infty(E^*)$
$$ \{\tilde\zvf,\tilde\zc\}_{\zL^{\widehat E^*}}=
{(\{\zvf,\zc\}_{\zL^{E^*}})}^\sim.
$$
\end{description}
\end{lemma}
Now, repeating the arguments from the classical case, one  easily  derives
the following.
\begin{theorem} For any bisection $J\in \A^2(E)$ of the vector bundle $E$
of a Jacobi algebroid $(E,\zF)$ the following  are equivalent:
\begin{description}
\item{(1)} $J$ is a canonical structure, i.e. $\lna  J,J\rnaf=0$;
\item{(2)} The canonical Jacobi bracket $J_\zF^{E^*}$
and $-\widehat J_\zF$ are $\sharp_J$-related;
\item{(3)}   The bivector fields  $\zL^{E^*}$   and   $-\widehat J_\zF^c$
are $\sharp_J$-related;
\item{(4)} For any `1-forms' $\zm,\zn\in\zW^1(E)$,
$$ \sharp_J([\zm,\zn]_J)=[\sharp_J(\zm),\sharp_J(\zn)],
$$
where the bracket on the right-hand-side is the Lie algebroid bracket on
$E$.
\end{description}
\end{theorem}
Note that a canonical structure for a Jacobi algebroid gives rise  to  a
{\it triangular Jacobi bialgebroid}  \cite{GM}  (or  a  {\it  triangular
generalized Lie bialgebroid} in the terminology of \cite{IM}).


\end{document}